\documentclass[11pt,a4paper]{article}
\usepackage[pagebackref=true,linktoc=page,colorlinks=true,linkcolor=blue,citecolor=blue]{hyperref}
\usepackage{setspace}
\usepackage{epsfig}
\usepackage{psfrag}
\usepackage[retainorgcmds]{IEEEtrantools} 
\usepackage{fullpage,color}
\usepackage{authblk}
\usepackage{latexsym}
\usepackage{amsmath,amssymb,amsfonts,amsthm}
\usepackage{graphicx}
\usepackage{natbib}
\usepackage{pstool,psfrag}
\usepackage[section]{placeins}

\newtheorem{lemma}{Lemma}

\newtheorem{corollary}{Corollary}
\newtheorem{proposition}{Proposition}
\newtheorem{remark}{Remark}
\makeatletter
\newcommand{\mbbR}{\mbox{$\mathbb{R}$}}

\definecolor{light-gray}{gray}{0.80} 
\title{\textbf{Conditioned limit laws for inverted max-stable
    processes}} \author[*]{\textbf{Ioannis Papastathopoulos}}
\author[$\dagger$]{\textbf{Jonathan A. Tawn}} \affil[*]{\small School
  of Mathematics and Maxwell Institute, University of Edinburgh,
  Edinburgh, EH9 3FD} \affil[$\dagger$]{\small Department of Mathematics
  and Statistics, Lancaster University, Lancaster, LA1 4YF}
\affil[$$]{\small i.papastathopoulos@ed.ac.uk $\quad$
  j.tawn@lancaster.ac.uk}  \date{}

\begin{document}

\maketitle
\begin{abstract}
  Max-stable processes are widely used to model spatial extremes.\
  These processes exhibit asymptotic dependence meaning that the large
  values of the process can occur simultaneously over space.\
  Recently,\ inverted max-stable processes have been proposed as an
  important new class for spatial extremes which are in the domain of
  attraction of a spatially independent max-stable process but instead
  they cover the broad class of asymptotic independence.\ To study the
  extreme values of such processes we use the conditioned approach to
  multivariate extremes that characterises the limiting distribution
  of appropriately normalised random vectors given that at least one
  of their components is large.\ The current statistical methods for
  the conditioned approach are based on a canonical parametric family
  of location and scale norming functions.\ We study broad classes of
  inverted max-stable processes containing processes linked to the
  widely studied max-stable models of Brown-Resnick, Schlather and
  Smith,\ and identify conditions for the normalisations to either
  belong to the canonical family or not.\ Despite such differences at
  an asymptotic level, we show that at practical levels, the canonical
  model can approximate well the true conditional distributions.
 \end{abstract}


\section{Introduction}
Extreme environmental events, such as hurricanes, heatwaves, flooding
and droughts, can cause havoc to the people affected and typically
result in large financial losses.\ The impact of this type of event is
often exacerbated by the event being severe over a large spatial
region.\ The statistical modelling of spatial extremes is a rapidly
evolving area \citep{davpadrib} and is crucial to understanding,
visualizing and predicting the extremes of stochastic processes.\ The
approach that is currently most used for modelling spatial extreme
values assumes the environmental process is a max-stable process
\citep{haan84}.\ The most widely used max-stable processes are the
Smith \citep{smit90b}, Brown-Resnick \citep{browresn77,kabletal09},
the extremal Gaussian \citep{schl02}, and the extremal-$t$
\citep{demarta2005,nikoetal09} processes.

Max-stable processes are the only non-trivial limit of point-wise
normalised maxima of independent and identically distributed
realisations of a stochastic processes.\ When max-stable processes are
observed at a finite number of locations their joint distribution is a
multivariate extreme value distribution,\ which is underpinned by the
assumption of the original variables satisfying the dependence
structure conditions of multivariate regular variation
\citep{resn87}.\ Max-stable processes have marginal generalised
extreme value distributions \citep{cole01} and a complex non-negative
dependence structure which has a restricted form.\ To understand this
restriction, let $\{X_M(s),s \in \mathbb{R}^2\}$ be a spatial
max-stable process with continuous marginal distribution functions
$F_{s}$ and corresponding inverse denoted by $F_s^{\leftarrow}$.\
Then, if for any $s_1,s_2\in \mbbR^2$ where $X_M(s_1)$ and $X_M(s_2)$
are not independent, it follows that the dependence coefficient
\[
  \chi=\lim_{p\rightarrow 1} \Pr\left\{X_M(s_2) > F_{s_2}^{\leftarrow}(p) \mid X_M(s_1)
    > F_{s_1}^{\leftarrow}(p) \right\},
\]
is positive.\ This property is termed asymptotic dependence;\ for
max-stable processes it implies that the spatial properties of extreme
events are independent of the severity of the event \citep{davpadrib}.

\cite{wadtawn12} introduced the class of inverted max-stable processes
and these are used in \cite{davhusthi13}.\ Any inverted max-stable
process $\{X(s),s\in \mathbb{R}^2\}$ with unit exponential margins,
i.e., for all $s\in \mathbb{R}^2$ and $x>0$, $\Pr(X(s) < x) =
1-\exp(-x)$, $x>0$, can be represented by
\[
  X(s)=1/X_F(s) \quad s\in\mathbb{R}^2,
\]
where $\{X_F(s), s\in \mathbb{R}^2\}$ is a max-stable process with
unit Fr\'{e}chet margins, i.e., for all $s\in\mathbb{R}^2$ and $x>0$,
$\Pr(X_F(s) < x) = \exp(-1/x)$.\ Thus, for all
$s_1,s_2\in\mathbb{R}^2$, the dependence structure between large
$X(s_1)$ and $X(s_2)$ is equivalent to the dependence structure
between small $X_F(s_1)$ and $X_F(s_2)$ and hence differs from the
max-stable form.\ All non-perfectly dependent inverted max-stable
processes are in the domain of attraction of spatially independent
max-stable processes (see Section \ref{sec:ms_and_ims}), meaning that
their point-wise normalised maxima are independent, i.e., for all
$s_1,s_2 \in \mathbb{R}^2$ , with $s_1\not= s_2$,
\begin{equation}
\lim_{n\rightarrow \infty} \Pr\left(\max_{i=1,\hdots,n}X_i(s_1) - \log
n <x_1,\max_{i=1,\hdots,n}X_i(s_2) - \log
n <x_2\right) = \prod_{i=1}^{2}\exp\{-\exp(-x_i) \},
\label{eq:componentmax}
\end{equation}
for any $x_1,x_2 \in \mathbb{R}$, where $\{X_i(s),s \in
\mathbb{R}^2\}$, $i=1,\hdots,n$ denotes a sequence of independent and
identically distributed inverted max-stable processes with unit
exponential margins.

To reveal the extremal dependence structure for asymptotically
independent random variables, alternative asymptotic properties need
to be studied.\ \cite{ledtawn96,ledtawn97} and \cite{resnick02}
explore the joint tail through the limiting joint survivor function,
\begin{equation}
  g(x_1,x_2)=\lim_{t\rightarrow \infty} \mathcal{L}(t)\exp\{t/\eta(s_1,s_2)\}
  \Pr\left\{X(s_1)>x_1+t, X(s_2)>x_2+t\right\}
\label{eq:eta}
\end{equation}
where $0<\eta<1$, known as the coefficient of tail dependence, and
$\mathcal{L}$ is a slowly varying function at $\infty$, are selected
so that $0<g(x_1,x_2)<\infty$,\ and the dependence structure revealed
by limit~\eqref{eq:eta} is known as hidden regular variation.\ A
weakness with this approach is that it fails to describe the behaviour
of the $X(s_2)$ values that occur with the largest values of
$X(s_1)$.\ Instead a conditioned approach is required which looks at a
more subtle normalisation for $X(s_2)$ that focuses on the region of
the joint distribution which is most likely when conditioning on
variable $X(s_1)$ being large.\ This is the approach we take in this
paper.


For a bivariate random variable $(X,Y)$ with unit exponential margins
and general dependence structure the conditioned extremes limit theory
of \cite{hefftawn04} and \cite{heffres07} is equivalent to the
assumption that there exist location and scaling norming functions
$a:\mathbb{R}_{+} \rightarrow \mathbb{R}$ and $b:\mathbb{R}_{+}
\rightarrow \mathbb{R}_{+}$, such that, for any $x>0$ and $z\in
\mbbR$,
\begin{equation}
  \lim_{u\rightarrow \infty}\Pr\left\{X-u>x,\{Y-a(X)\}/{b(X)}<z\mid
    X>u\right\} = \exp(-x)G(z),
  \label{eq:cond_rep}
\end{equation}
where $G$ is a non-degenerate distribution function.\ To ensure $a$,\
$b$ and $G$ are uniquely defined the condition $\lim_{z\rightarrow
  \infty}G(z) = 1$ is required,\ so $G$ places no mass at $+\infty$
but some mass is allowed at $-\infty$.\ For positively dependent
random variables, \cite{hefftawn04} found that, for all the standard
copula models studied by \cite{joe97} and \cite{Nelsen06}, the norming
functions $a(x)$ and $b(x)$, fell into the simple canonical parametric
family
\begin{equation}
a(x) = \alpha  x \quad
\text{and} \quad b(x) =
x^{\beta},
\label{eq:HT_class}
\end{equation}
where $\alpha \in [0,1]$ and $\beta\in (-\infty,1)$.\ The case $\alpha
= 1$ and $\beta = 0$ corresponds to $\chi>0$, whereas any other
combination of $\alpha$ and $\beta$ gives $\chi=0$.\ With standard
Laplace margins, \cite{keefpaptawn13} extended model
class~(\ref{eq:HT_class}) to $-1\leq \alpha \leq 1$ to account for
negatively dependent random variables when $\alpha \in[-1,0)$.\
However,\ this is not relevant for inverted max-stable processes as
they are non-negatively dependent.

A key question is why have we made the restriction of unit exponential
margins on the inverted max-stable process.\ In the style of copula
methods \citep{Nelsen06} we assume identical margins.\ Conditioned
limit theory has studied limiting presentation \eqref{eq:cond_rep}
with the margins $(X,Y)$ taken as identically distributed Gumbel
variables \citep{hefftawn04},\ which is asymptotically equivalent to
unit exponential margins but mathematically less clean,\ whereas
\citet{eastoetawn12} take identically distributed generalised Pareto
distributions.\ In contrast \cite{heffres07} work with $X$ in the
domain of attraction of the generalised extreme value distribution
distribution, but do not impose any constraint on $Y$, and find that
it is not always possible to achieve an affine normalisation as in
limit~\eqref{eq:cond_rep} after marginal transformation to identical
margins.\ \cite{kulisoul14} consider $(X,Y)$ with margins that have
regularly varying tails,\ with location function $a(x)=0$.\ Through
the paper we also find that for broad and important classes of
inverted max-stable processes, it is possible to achieve affine
normalisations after transformation to identical exponential marginal
variables.\ Furthermore, for some of these classes of inverted
max-stable distributions limit~\eqref{eq:cond_rep} does not hold with
margins that have regularly varying tails but holds with unit exponential margins.\
Thus our restriction to unit exponential margins provides all the
necessary ground work for deriving conditioned limits with any
marginal distributions for inverted max-stable processes.\ In
Section~\ref{sec:exp_scale} we discuss the implications of the
marginal choice, and in particular show how once the limit
relationship has been derived for unit exponential margins alternative
limit results follow immediately for other marginal choices.\

Canonical family~(\ref{eq:HT_class}) has been subject of criticism
\citep{Smit04} since the functions $a$ and $b$ seem to be `proof by
example' rather than a general result.\ \cite{heffres07} note that
under assumption (\ref{eq:cond_rep}) and in unit exponential
marginals,
\[
  \lim_{t\rightarrow \infty}b(t+x)/b(t) = \psi_1(x) \quad
  \text{and}\quad \lim_{t\rightarrow \infty}\left\{a(t + x) -
    a(t)\right\}/b(t) = \psi_2(x),
\]
for any $x\in\mathbb{R}$, where $\psi_1$ and $\psi_2$ are real
functions.\ For the canonical family~(\ref{eq:HT_class}), it readily
follows that $(\psi_1,\psi_2)= (1,0)$ if $\beta > 0$ and
$(\psi_1,\psi_2)=(1,\alpha x)$ if $\beta = 0$ and this condition is
also satisfied by a range of regularly varying functions.\ However, no
examples have been published to date other than the canonical
form~(\ref{eq:HT_class}).

For simplicity we focus on bivariate characterisations of the inverted
max-stable process,\ with all joint distributions following inverted
bivariate extreme value distributions.\ For the rest of the paper we
denote by $(X,Y)$ a bivariate random variable with inverted bivariate
extreme value distribution with unit exponential margins and derive
$a$,\ $b$ and $G$.\ We find classes of this family where the
normalisation required to achieve property~(\ref{eq:cond_rep}) either
falls in the canonical family~(\ref{eq:HT_class}) or a more general
form is required.\ Examples of the former include the inverted
max-stable models with Schlather \citep{schl02} and extremal-$t$
\citep{demarta2005,nikoetal09} dependence models, and the latter
include inverted max-stable model with Smith \citep{smit90b}
dependence model.\ We show that the distinction between these classes
is determined by the behaviour of the spectral measure of the
underlying max-stable process near its lower end point.

The statistical conditioned model of \cite{hefftawn04} assumes that
the limiting relationship~(\ref{eq:cond_rep}) holds exactly for all
values $X>u$ for a suitably high threshold $u$ with the functions $a$
and $b$ in canonical form~(\ref{eq:HT_class}).\ This model has been
found to fit well in various applications
\citep{pauletal06,keefST09,hila11,eastoetawn12,papatawn14}.\ Our
identification of the existence of the new classes that do not fall in
canonical family~(\ref{eq:HT_class}) questions the validity of the
generic use of the canonical family for statistical modelling.\
Therefore,\ we compare the new models with the current statistical
approach of \cite{hefftawn04} and show,\ through simulation,\ that at
practical levels,\ that the use of canonical family gives good
approximations to the conditional distribution of the inverted
max-stable model with Smith dependence and highlight examples where a
good approximation does not hold.

The paper is structured as follows.\ In
Section~\ref{sec:inverted_models} we present the classes of max-stable
and inverted max-stable distributions and explore the implications of
deriving results on general margins.\ In
Section~\ref{sec:new_examples} we present the conditional
representation of the class of inverted max-stable distributions with
spectral densities of the associated max-stable process being
regularly varying and $\Gamma$-varying spectral densities at their
lower end-point.\ In Section~\ref{sec:spatial} we discuss the
spatial extension of our results.\ Our derivations and proofs are
included in the Appendix.

\section{Bivariate inverted max-stable distributions}
\label{sec:inverted_models}
\subsection{Max-stable and inverted max-stable distributions}
\label{sec:ms_and_ims}
Max-stable distributions arise naturally as the only non-degenerate
limit distributions of appropriately normalised component-wise maxima
of random vectors.\ In unit Fr\'{e}chet margins, and for $x,y > 0$, a
bivariate max-stable distribution function is defined by
\begin{equation}
  F(x,y)=\exp\left\{-V\left(x,y\right)\right\} =
  \exp\left[- \int_{0}^{1}\max\left\{w/x,(1-w)/y\right\}\,\mathrm{d}H(w)\right],
\label{eq:max_stable}
\end{equation}
where $V$ is termed the exponent measure and $H$ is an arbitrary
finite measure on $[0,1]$, known as the spectral measure, with total
mass 2, satisfying the marginal moment constraint $\int_0^1 w\,\mathrm{d}H(w)=1$.\
\cite{coletawn91} showed that if $F$ has a density, then $H(w)$ has
spectral density $h(w)$ on the interior $(0,1)$ and can have mass
$H(\{k\})$, $k=0,1$, on each of $\{0\}$ and $\{1\}$, given by
\begin{IEEEeqnarray*}{rCl}
  h(w) &=& -\frac{\partial^{2}V}{\partial x\partial y}(w,1-w) \quad
  0<w<1,\nonumber
  \\  \\
  H(\{0\}) &=& -y^2\lim_{x \rightarrow 0} \frac{\partial V}{\partial
    y}(x,y), \quad \text{and} \quad H(\{1\}) = -x^2\lim_{y \rightarrow 0}
  \frac{\partial V}{\partial x}(x,y).\nonumber
\end{IEEEeqnarray*}
As the class of bivariate max-stable distributions does not admit a
finite dimensional parameterisation, a natural method for modelling
the spectral measure $H$ of expression~(\ref{eq:max_stable}) relies on
constructing parametric sub-classes of models that are flexible enough
to approximate any member from the class
\citep{coletawn91,ballschla11}.\ Two such sub-models are
\cite{huslreis89} and \cite{schl02} max-stable distributions which
have exponent measures, for $x,y>0$,
\begin{IEEEeqnarray}{rCl}
  &V(x,y)=\frac{1}{x}\Phi\left\{ \frac{\lambda}{2} +
    \frac{1}{\lambda}\log\bigg(\frac{y}{x}\bigg) \right\} +
  \frac{1}{y}\Phi\left\{ \frac{\lambda}{2} + \frac{1}{\lambda}
    \log\left(\frac{x}{y}\right) \right\} \quad \lambda\in (0,\infty), \label{eq:hr_model}\\\nonumber\\
  & V(x,y)=\frac{1}{2}\left(\frac{1}{x} + \frac{1}{y}\right) \left[1 +
    \left\{1 - 2 \left(1+\rho\right)\frac{x
        y}{\left(x+y\right)^2}\right\}^{1/2}\right] \quad \rho \in
  (-1,1) ,&\label{eq:schl_model}
\end{IEEEeqnarray}
respectively, where $\Phi$ is the cumulative distribution function of
the standard normal distribution.\ These are the exponent measures of
the pairwise distributions for the Smith and Schlather max-stable
models respectively.\ The parameters $\lambda$ and $\rho$ control the
strength of dependence.\ In particular, increasing and decreasing
values of $\rho$ and $\lambda$, respectively, imply stronger
dependence between $X$ and $Y$.

Given a max-stable distribution with exponent measure $V$ as in
equation~(\ref{eq:max_stable}), the bivariate random variable $(X,Y)$
follows the inverted max-stable distribution with unit exponential
margins if, for $x,y>0$, its joint survivor function is,
\begin{equation}
  \Pr\left(X>x,Y>y\right) =
  \exp\left\{-V\left(1/x,1/y
    \right)\right\}.
\label{eq:inv_max_stable}
\end{equation}
As 
\[
\Pr\left(X>x,Y>x\right)  =\exp(-x/\eta)= \{\Pr(X>x)\}^{1/\eta}, \quad
\eta=1/V(1,1) \in [1/2,1],
\]
the inverted max-stable distributions are either perfectly dependent
when $V(1,1)=1$ or asymptotically independent when $V(1,1)>1$.\
Specifically $\eta$ is the coefficient of tail dependence
\citep{ledtawn97}.\ This property explains the independence in
limit~\eqref{eq:componentmax}.\

A range of results are available to study the conditional
limit~\eqref{eq:cond_rep}.\ In particular, \cite{heffres07},
\cite{resnzebe14} and \cite{wadetal14} show that under various
conditions, all of the followings limits are identical to $G(z)$
in limit~\eqref{eq:cond_rep}:
\begin{align}
  \lim_{u\rightarrow \infty}\Pr\left\{ \frac{Y-a(u)}{b(u)}<z \big | X
    >u\right\} = \lim_{u\rightarrow \infty}\Pr\left\{
    \frac{Y-a(X)}{b(X)}<z \big | X >u\right\} = \lim_{u\rightarrow
    \infty}\Pr\left\{ \frac{Y-a(u)}{b(u)}<z \big | X = u\right\}.
\label{eq:reg_cond}
\end{align}
These conditions hold for all inverted max-stable distributions if the
associated spectral measure of the max-stable process places no point
mass on the interval $(0,1)$,\ i.e.,\ when $(X,Y)$ have a joint
density.\ So we can use any of these limits to derive the forms of $a,
b$ and $G$.\ Motivated by statistical considerations, we find that the
use the last expression is most simple to use.\ However,\ this is the
most restrictive in general as it requires the assumption of a joint
density.\ As all the parametric max-stable models have joint densities
so do the associated inverted max-stable distributions, so for us this
is not restrictive.\ For using this third limit form it is helpful to
note that the conditional survivor function is
\begin{equation}
  \Pr\left(Y>y\mid X=x\right) = 
  -V_{1}\left(1,x/y\right) \exp\left\{x - x
    V\left(1,x/y\right)\right\}  \quad y > 0,
\label{eq:approx_surv}
\end{equation}
where $V_{1}(x,y)=\partial V(x,y)/\partial x$.

\subsection{Conditional representation with different marginal distributions}
\label{sec:exp_scale}
Let $(X,Y)$ be a bivariate random variable with common unit
exponential margins and assume that limit~\eqref{eq:cond_rep} holds.\
Here we consider what this representation then implies for the
extremal conditioned distribution of $Y_{H_2} \mid X_{H_1}$, where the
bivariate random variable $(X_{H_1},Y_{H_2})$ has continuous marginal
distribution functions $H_1$ and $H_2$ respectively and has identical
copula to $(X,Y)$.\ For $i=1,2$,\ let $K_i(y) = -\log\{1-H_i(y)\}$ for
$y\in \mathbb{R}$ and denote its inverse by
$K_i^{\leftarrow}(y)=H_i^{\leftarrow}\{1-\exp(-y)\}$ for $y>0$.\ Then
$(X_{H_1},Y_{H_2}) = \{K_1^{\leftarrow}(X),K_2^{\leftarrow}(Y)\} $
have the required joint distribution and from
limit~\eqref{eq:cond_rep} there exist functions
$a:\mathbb{R}_+\rightarrow \mathbb{R}$ and $b:\mathbb{R}_+\rightarrow
\mathbb{R}_+$ such that for all $z\in \mbbR$
\begin{equation}
  \lim_{u\rightarrow \infty} \Pr\left(Y_{H_2} <
    K_2^{\leftarrow}\left[a\{K_1(X_{H_1})\} + b\{K_1(X_{H_1})\} z\right]\mid X_{H_1} >
    u\right)=G(z)
\label{eq:different}
\end{equation}
for non-degenerate $G(z)$.\ Therefore,\ if we can find a
location-scale normalisation when working with identical unit
exponential margins then limit~\eqref{eq:different} shows that these
results directly provide the appropriate conditioned limit for
non-identically distributed marginals.\ Furthermore,\
limit~\eqref{eq:different} shows that in general margins a
location-scale normalisation is not always possible even when it can
be achieved with unit exponential margins.\ Of course the converse is
true, but as we will see in Section~\ref{sec:new_examples} for the
class of inverted max-stable distributions limit~\eqref{eq:cond_rep}
holds with unit exponential margins, so limit~\eqref{eq:different} is
useful to give the conditioned limits in other marginals.

To help understand the implications of limit~\eqref{eq:different} it
is helpful to focus on specific forms for on $H_1$ and $H_2$.\
\citet{eastoetawn12} present limit~\eqref{eq:different} with
identically distributed generalised Pareto distributions.\
\cite{kulisoul14} work with regularly varying tails.\ Focusing on the
specific case of Pareto margins with $H_i(y)=1-y^{-\alpha_i}$ for
$\alpha_i>0,\ y>1$ and $i=1,2$, limit~\eqref{eq:different} then
becomes
\[
  \lim_{u\rightarrow \infty} \Pr\left\{Y_{H_2} <
    \exp\left(\left[a\{\alpha_1\log(X_{H_1})\} 
        + b\{\alpha_1\log(X_{H_1})\} z\right]/{\alpha_2}\right) \,\Big|\, X_{H_1} 
    > u \right\}=G(z)
\]
so a location-scale normalisation in these margins can only be
achieved with $b(y)\equiv 1$ and then only a scaling is required.\ Thus
studying limit~\eqref{eq:cond_rep} with regularly varying tails only
requires a scaling but critically cannot cover any cases where the $b$
scaling function when using unit exponential margins differs from
$b(y)\equiv 1$.\ As we will see in Section~\ref{sec:new_examples}, for the
class of inverted max-stable distributions with unit exponential
margins we have many classes where $b(y)\not= 1$.\ Therefore, to reveal
the full structure of the conditioned extremal behaviour of inverted
max-stable distributions revealed by limit~\eqref{eq:cond_rep}, it is
essential to work in marginal variables which are tail equivalent to
the unit exponential.\ Hence for the remainder of the paper we work
exclusively with unit exponential margins acknowledging that these
results apply directly with different margins using
limit~\eqref{eq:different}.

\section{Conditional representations}
\label{sec:new_examples}

\subsection{Known representations}
\label{sec:known}
\cite{hefftawn04} explored the conditional
representation~(\ref{eq:cond_rep}) for the class of inverted
max-stable distributions subject to the assumption that $H$ places all
the mass in $(0,1]$ and that the spectral density is regularly varying
\footnote{A function $h:\mbbR_+\rightarrow \mbbR_+$ is regularly
  varying at $0$, with index $t\in \mbbR$, short-hand $h\in
  R_{t}(0^+)$ if, for all $w>0$, $ \lim_{s\rightarrow 0^+} h(s w)/h(s)
  = w^{t}$.\ For any $h\in R_{t}(0^+)$, it follows that for all $w>0$,
  $h(w)=w^t\mathcal{L}(w)$, where $\mathcal{L}$ is a slowly varying
  function, i.e.,\ $\mathcal{L}\in R_0(0^+)$.} 
\begin{IEEEeqnarray}{rCl}
  &h(w)\sim \mathcal{L}(w-w_{\ell}) (w-w_{\ell})^{t} \quad \text{as
    $w\rightarrow w_{\ell}\in[0,1/2)$},&\label{eq:h_HT}
\end{IEEEeqnarray}
for $w_{\ell} = 0$ and $\mathcal{L}(w)$ slowly varying at 0 with
$\lim_{w\rightarrow w_{\ell}^+}\mathcal{L}(w-w_{\ell}) = s > 0$ and $t>-1$.\
Later, we consider more general formulations with
$\mathcal{L}(w)\rightarrow \infty$ as $w\rightarrow 0^+$ and $w_{\ell} >
0$.\ Under setting~(\ref{eq:h_HT}), the
normalisation~(\ref{eq:HT_class}) required to give a non-degenerate
limiting conditional law has $\alpha=0$, $\beta=(t+1)/(t+2)$ and the
limit is of Weibull type, i.e.,\
\begin{equation}
  \lim_{u\rightarrow \infty}\Pr\left(Y<X^{\beta}z\mid X>u\right) = 1-
  \exp\left\{-\frac{s z^{t+2}}{(t+1)(t+2)}\right\}, \quad \text{for $z>0$.}
\label{eq:reg_var_lim}
\end{equation}
Examples of inverted max-stable models
satisfying~\eqref{eq:reg_var_lim} include those with logistic and
Dirichlet dependence structure \cite{coletawn91}.

\subsection{Regular variation at lower tail of spectral measure}
\label{sec:reg_var}
We derive the conditional representation~(\ref{eq:cond_rep}) for the
class of inverted max-stable distributions covering more general
spectral measures than those studied in Section~\ref{sec:known}.\ In
particular,\ we explore model~(\ref{eq:h_HT}) and the effect on the
normalising functions $a(x)$ and $b(x)$ when the spectral measure $H$
places its mass in a sub-region,\ $[w_{\ell},w_{u}]$ say,\ of $[0,1]$.\
Motivated by the Schlather distribution~(\ref{eq:schl_model}),\ for
which the spectral measure places mass at $\{0\}$,\ we also explore
the assumption of possible mass on the lower end point $w_{\ell}$ of $H$,\
but no point mass at any other $w\in(w_{\ell},w_{u})$.\ Although having mass
at $w_{\ell}$, if $w_{\ell}>0$, implies that $(X,Y)$ do not have a joint
density, as there is a singular component on the boundary of the
sample space of $(X,Y)$, if this boundary is avoided
then~\eqref{eq:approx_surv} is still valid.\ The results in the rest
of the paper even hold if there is mass at some point in $(w_{\ell},w_{u})$,
with modified proofs, but we avoid this unnecessary generalisation.
\begin{lemma}
  \label{le:derivative}
  Let $w_{\ell}$,\ $w_{u}$, be the lower and upper end
  points,\ respectively,\ of the spectral measure $H$ of an inverted
  max-stable distribution~(\ref{eq:inv_max_stable}),\ i.e.,\
  \[
  w_{\ell} = \inf \left\{ {0 \leq w < 1/2}: H(w)\geq0\right\},\quad w_{u} =
  \sup\left\{ 1/2 < w \leq 1: H(w) \leq 2\right\},
  \]
  and assume that, apart from the points $w_{\ell}$ and $w_{u}$, for which
  $H(w_{\ell}),H(w_{u})\geq 0$, 
  the spectral measure is absolutely continuous with respect to
  the Lebesgue measure.\ Then,\ if $y/(x + y) \downarrow w_{\ell}$,\ for
  $y=y(x)$ a function of $x$, as $x\rightarrow \infty$,
  \begin{IEEEeqnarray}{rCl}
    \label{eq:partial_limit}
    V_{1}\left(1,x/y\right) &\rightarrow& w_{\ell} H(\{w_{\ell}\}) - 1.
  \end{IEEEeqnarray}
\end{lemma}
In Lemma~\ref{le:derivative} the case of perfect positive dependence
between $X$ and $Y$,\ i.e.,\ $w_{\ell} = w_{u} = 1/2$,\ is excluded since
there can be no possible normalisation such that $G$ in
expression~(\ref{eq:cond_rep}) is non-degenerate.\
Proposition~\ref{prop:log_survivor_cases} gives the asymptotic form of
the log-conditional survivor function of $Y\mid X=x$, for large $x$.
\begin{proposition}
  \label{prop:log_survivor_cases}
  Under the conditions of Lemma~\ref{le:derivative} and for $h(w)$ as
  in expression~(\ref{eq:h_HT}), for $y/(x+y)\downarrow w_{\ell}\in
  [0,1/2)$ for $y=y(x)$ a function of $x$, as $x \rightarrow \infty$,\
  we obtain that $\log\Pr\left(Y>y \mid X=x\right)$ is asymptotically
  equivalent to

\noindent $(i)$ for $H(\{w_{\ell}\})=0$,
\begin{IEEEeqnarray}{rCl}
 && - x \mathcal{L}\left(\frac{y}{x + y} - w_{\ell}\right)
  \left(\frac{y}{x + y} - w_{\ell}\right)^{t+2}\big/\{(1-w_{\ell})(t+1)(t+2)\},
\label{eq:nomass_cond_distr}
\end{IEEEeqnarray}
 
\noindent $(ii)$ for $H(\{w_{\ell}\}) > 0$, 
\begin{IEEEeqnarray}{rCl}
  && \log\left\{1 - w_{\ell} H(\{w_{\ell}\})\right\} - (x+y)\left(\frac{y}{x+y} -
    w_{\ell}
  \right) H(\{w_{\ell}\}).
\label{eq:mass_cond_distr}
\end{IEEEeqnarray}
\end{proposition}


As there is no contribution from the spectral density in
expression~(\ref{eq:mass_cond_distr}),\ a general form for the
normalisation can be obtained directly.\ General forms of the
normalising functions $a(x)$ and $b(x)$ cannot be obtained from
representation~(\ref{eq:nomass_cond_distr}) without additional
assumptions so it is helpful to consider a condition on
the slowly varying function $\mathcal{L}$ at $0$.\ \cite{hefftawn04}
results were based on the case of $\mathcal{L}$ having a finite right
limit at $0$,\ Corollary~\ref{cor:normings_reg_var} extends this.
\begin{corollary}
\label{cor:normings_reg_var}
Under the conditions of Proposition~{\ref{prop:log_survivor_cases}},\
limit~\eqref{eq:cond_rep} holds as follows

\noindent $(i)$ for $H(\{w_{\ell}\})=0$, then $a(x)=\{w_{\ell} /(1-w_{\ell})\}x$ and
$b(x)=x^{(t+1)/(t+2)}\mathcal{L}\left\{x^{-1/(t+2)}\right\}^{-1/(t+2)}$
for $x>0$ and assuming that
\begin{equation}
  \lim_{w \rightarrow 0^+} \frac{\mathcal{L}\left\{w
      \mathcal{L}(w)^{-\tau}\right\}}{\mathcal{L}(w)} = 1,
\label{eq:slow_condition}
\end{equation}
for all $\tau \in (0,1)$, then
\[
G(z) = 1 -
\exp\left\{-\frac{ (1-w_{\ell})^{3+2t} z^{t+2}}{(t+1)(t+2)} \right\} \quad
\text{for $z>0$}.
\]

\noindent $(ii)$ For $H(\{w_{\ell}\})>0$, then $a(x)=\{w_{\ell}/(1-w_{\ell})\}x$ and
$b(x)=1$ for $x>0$, then 
\[
G(z) = 1 - \left\{1- w_{\ell}H(\{w_{\ell}\})\right\}\exp\left\{-(1-w_{\ell})
  {H(\{w_{\ell}\})} z\right\} \quad \text{for $z > 0$},
\]
and $G(z)=0$, for $z < 0$, so $G(\{0\})= w_{\ell} H(\{w_{\ell}\})$.
\end{corollary}

Condition~(\ref{eq:slow_condition}) is satisfied by a range of slowly
varying functions,\ including those studied in Section 3.1 as
well as by functions that approach $\infty$ when the argument tends to
zero.\ Examples for $\mathcal{L}(w)$,\ with $\lim_{w \rightarrow
  0^+}\mathcal{L}(w) = \infty$,\ satisfying
condition~(\ref{eq:slow_condition}) include $\log_{\kappa}\left(-\log
  w\right)$,\ $\kappa \in \mathbb{N}_0$,\ $ \exp\left\{(-\log
  w)^{\nu}\right\}$,\ $\nu \in (0,1/2)$ and $ \exp\left\{ - \log w /
  \log\left(-\log w\right)\right\}$,\ where $\log_k$ is the iterated
logarithm function defined recursively by $\log_k x =
\log_{\kappa-1}\log x$,\ $\log_0 x = x$ and $\log_1 x =\log x$.
\begin{remark}
  For $\lim_{w\rightarrow 0^+}\mathcal{L}(w) = s > 0$,\ all cases of
  norming functions $(i)$-$(ii)$ in
  Corollary~\ref{cor:normings_reg_var} reduce,\ after absorbing $s$
  into the limiting law,\ to the parametric class of Heffernan--Tawn,\
  i.e.,\ $a(x)+b(x)z=\alpha x + x^{\beta} z$,\ where
  $\alpha=w_{\ell}/(1-w_{\ell}) \in [0,1)$ for $w_{\ell}\in[0,1/2)$ and
  $\beta=(t+1)/(t+2)\in [0,1)$ for $t\geq -1$.\ This is the first
  example to be known with both $\alpha$ and $\beta$ a function of
  parameters.\ 
\end{remark}
\begin{remark}
When $\lim_{w\rightarrow 0^+}\mathcal{L}(w) = \infty$ and subject to
condition~(\ref{eq:slow_condition}),\ the additional factor
$\mathcal{L}\{x^{-1/(t+2)}\}^{-1/(t+2)}$ enters in the scaling
function and reduces the rate of increase of $b(x)$ to $\infty$,\ as
$x\rightarrow \infty$.
\end{remark}

Another interesting case which is satisfied by many max-stable models
that appear in the literature,\ is when $w_{\ell}=0$ and $H(\{w_{\ell}\})>0$ in
case $(ii)$ of Corollary~\ref{cor:normings_reg_var}.\ In this case,\
$a(x)=0$ and $b(x)=1$ for all $x>0$ and the random variables $X$ and
$Y$,\ conditionally on $X>u$,\ are near independent in the terminology
of \cite{ledtawn97},\ with exact independence occurring when the
limit distribution is unit exponential,\ i.e.,\ when $H(\{w_{\ell}\})=1$.\
Two such max-stable models come from the extremal-$t$
\citep{nikoetal09} and Gaussian-Gaussian \citep{wadtawn12} processes,\
for which the exponent measures of their bivariate distributions are
\begin{IEEEeqnarray}{rCl}
  V(x,y)& = &
  \frac{1}{x}T_{\nu+1}\left[\frac{\left(y/x\right)^{1/\nu}-\rho}{\left\{
        (1-\rho^2)/(\nu + 1)\right\}^{1/2}}\right]+   \frac{1}{y}T_{\nu+1}\left[\frac{\left(x/y\right)^{1/\nu}-\rho}{\left\{
        (1-\rho^2)/(\nu + 1)\right\}^{1/2}}\right],\label{eq:extremal_t}
  \\\nonumber \\
  V(x,y)&=&\frac{1}{2}\left(\frac{1}{x}+\frac{1}{y}\right) +
  \frac{1}{2}\int_{\mbbR^2}\left\{\frac{\phi_2\left(u\right)^2}{x^2} -2
    \rho(h)\frac{\phi_2\left(u\right) \phi_2\left(h - u\right)}{x y} +
    \frac{\phi_2\left(h - u\right)^2}{y^2}\right\}^{1/2}du,
  \label{eq:Gaussian_Gaussian}
\end{IEEEeqnarray}
respectively,\ where $\nu>0$, $h\in \mbbR^2_+$,\ $\rho(h)\in [-1,1]
$ is a valid correlation function,\ $T_{\nu+1}$ is the distribution
function of the standard-$t$ distribution with $\nu+1$ degrees of
freedom,\ and $\phi_2$ is the density of the standard bivariate normal
distribution with correlation $\rho(h)$.\ The corresponding mass on the
lower end point $w_{\ell}=0$ of models~(\ref{eq:extremal_t})
and~(\ref{eq:Gaussian_Gaussian}) is
\[
H(\{0\}) = T_{\nu+1}\left[-\rho\left(\frac{\nu + 1}{1 -
      \rho^2}\right)^{1/2}\right] \quad \text{and} \quad
H(\{0\})=\frac{1-\rho(h)}{2},
\]
respectively.\ Table~\ref{tab:mass} gives a collection of other
max-stable models,\ including the Schlather
distribution~(\ref{eq:schl_model}),\ placing positive mass on $\{0\}$.

\begin{table}[htbp!]
  \centering
  \caption{The mass of the spectral measure
    on $\{0\}$ of bivariate exponent measures,\ from top to bottom,\ of mixed,\ 
    asymmetric and asymmetric mixed logistic \citep{tawn88},\ Schlather \citep{schl02}
    and \cite{marsolki67a} distributions.\ The final column shows the parameter space,\ $\Theta$,
    of the model.}
  \begin{tabular}{c c c}
    \hline
    $V(x,y)$& $H(\{0\})$ & $\Theta$\\\hline
    $\left(\frac{1}{x}+\frac{1}{y}\right) - \frac{\theta}{x+y}$& $1-\theta$ & $\theta \in (0,1)$\\ 
    $\frac{1-\theta}{x} + \frac{1-\phi}{y} + \left\{\left(\theta/x\right)^{1/\alpha} + \left(\phi/y\right)^{1/\alpha}\right\}^\alpha$& $1-\phi$ & $0 \leq \theta, \phi, \alpha \leq 1$ \\  
    $\frac{1}{x}+\frac{1}{y} -\frac{1}{x y}\left(\frac{1}{x}+\frac{1}{y}\right)^{-2}\left(\frac{\theta + \phi}{x} + \frac{2\phi + \theta}{y}\right)$ & $1- \phi - \theta$ & $\theta,\theta+3\phi>0$ and $\theta + \phi, \theta + 2\phi \leq 1$  \\ 
    $\frac{1}{2}\left(\frac{1}{x} + \frac{1}{y}\right) \left[1 +
      \left\{1 - \frac{2 \left(1+\rho\right) x
          y}{\left(x+y\right)^2}\right\}^{1/2}\right] $& $\left(1-\rho\right)/2$ & $\rho \in (-1,1)$\\ 
    $\alpha \left(\frac{1}{x}+\frac{1}{y}\right)+(1-\alpha)\max\left\{1/x,1/y\right\}$& $\alpha$ & $0\leq\alpha\leq 1$ 
  \end{tabular}
  \label{tab:mass}
\end{table}
\subsection{Inverted Smith model}
\label{sec:inv_HR}
In this section we focus on the limiting conditional representation of
the inverted Smith model~(\ref{eq:hr_model}).\ This model has the same
bivariate copula as the H\'{u}sler--Reiss distribution.\ The spectral
measure of the Smith max-stable distribution places no mass on any
$0\leq w \leq 1$,\ and the spectral density satisfies
\begin{IEEEeqnarray}{rCl}
  &h(w) \sim \frac{\exp(-\lambda^2/8)}{\lambda\, (2\pi)^{1/2}}
  w^{-3/2}\exp\left\{-\left(\log w\right)^2/(2\lambda^2)\right\}\quad
  \text{as $w\rightarrow 0$}.  &\label{eq:h_hr}
\end{IEEEeqnarray}
This corresponds to a different form than expression~(\ref{eq:h_HT})
or its more general forms of the slowly varying function
$\mathcal{L}$.\ In particular, the spectral density is
$\Gamma$-varying\footnote{A function $h$ is $\Gamma$-varying at $0$
  \citep{haan70} with auxiliary function $f$, short-hand $h\in
  \Gamma_{f}(0^+)$, if for all $w\in \mathbb{R}$, $ \lim_{s\rightarrow
    0^+} h\left\{s + w f(s)\right\}/h(s) = \exp(w)$.} at 0 with
auxiliary function
\[
f(w) = -\lambda^2 w/\log w.
\] 
As Proposition~\ref{prop:HR} shows,\ this example leads to a different
form for the normalising functions $a(x)$ and $b(x)$ than the ones
considered by \cite{hefftawn04}.
\begin{proposition}
  \label{prop:HR}
  Assume that $(X,Y)$ follows the inverted max-stable
  distribution~(\ref{eq:inv_max_stable}) with exponent
  measure~(\ref{eq:hr_model}).\ Limit~\eqref{eq:cond_rep} holds with,
  for $x>0$ and $z\in\mathbb{R}$,
  \begin{equation}
    a(x)= x\exp\left\{- \lambda \,(2\log x)^{1/2}
      + \frac{\lambda \log \log x}{(2 \log x)^{1/2}} +  \frac{\lambda^2}{2}\right\},\, b(x) = a(x)/(\log x)^{1/2}\label{eq:a_hr},
  \end{equation}
  and 
  \begin{equation}
    G(z) = 1 - \exp\left\{ -
      \frac{\lambda}{(8\pi)^{1/2}}\exp\left(\sqrt{2}z/\lambda\right)\right\}.
\label{eq:limit_hr}
  \end{equation}
\end{proposition}
\begin{remark}
  Limit distribution~(\ref{eq:limit_hr}) is of reverted Gumbel type,\
  which is different from the limits in
  Corollary~\ref{cor:normings_reg_var}.\ The rate of convergence to
  the limit is order $\log\log u/(\log u)^{1/2}$.\
\end{remark}

\begin{figure}[htpb!]
  \centering
   \includegraphics[scale=1]{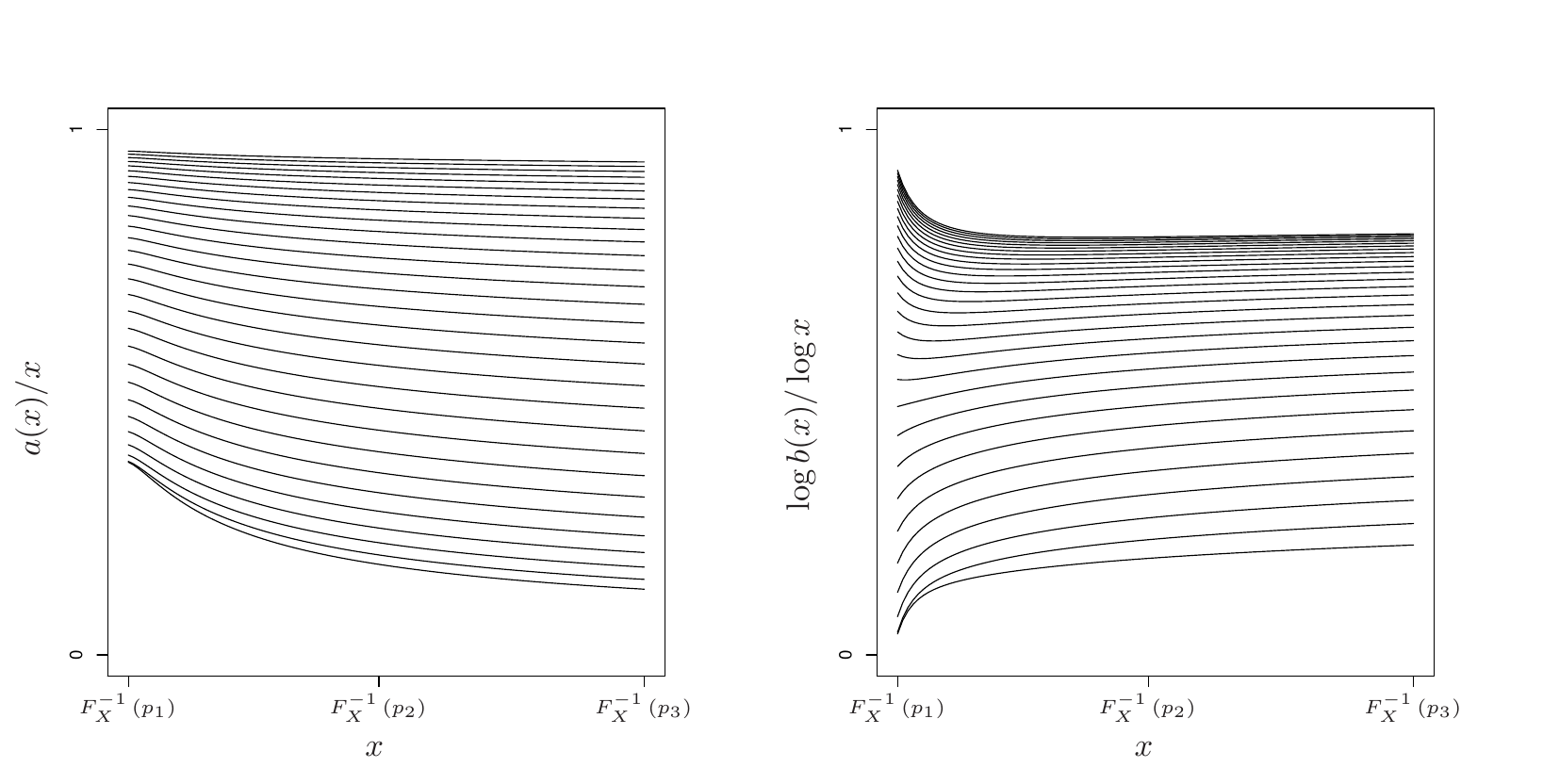}
   \caption{Plots of $a(x)/x$ and $\log b(x)/\log x$,\
     $x>F_X^{-1}(0.87)$,\ where $a(x)$ and $b(x)$ are given by
     expression~(\ref{eq:a_hr}),\ for different values of $\lambda$
     ranging from $0.01$ (bottom curves) to $20$ (top curves).\ The
     inverse of the unit exponential distribution function is
     $F_X^{\leftarrow}(p)=-\log(1-p)$,\ for $p\in(0,1)$.\ The values of
     $p_1$,\ $p_2$ and $p_3$ are 0.95,\ $1-10^{-7}$ and $1-10^{-13}$,
     respectively.}
\label{fig:norm_functions}
\end{figure}

A natural question that arises from this counter-example relates to
how well can the conditioned dependence model of \cite{hefftawn04}
with canonical family~(\ref{eq:HT_class}) approximate the conditional
distribution of $Y\mid X>u$,\ for large $u$,\ when the random vector
$(X,Y)$ follows the inverted max-stable distribution with Smith
dependence structure and unit exponential margins.\ To facilitate
comparisons between the two models,\ Figure~\ref{fig:norm_functions}
shows the graphs of $a(x)/x$ and $\log b(x)/\log x$ where $a(x)$ and
$b(x)$ are given by expression~(\ref{eq:a_hr}),\ for several values of
the dependence parameter $\lambda$ and a range of $x$ values above the
0.87 unit exponential quantile.\ Both plots show that $a(x)/x$ and
$\log b(x)/\log x$ are approximately constant for large $x$ so that
the canonical class of norming functions is likely to approximate well
$a(x)$ and $b(x)$ by $\alpha x$ and $x^\beta$,\ respectively.\
Subsequently,\ we simulated data from the inverted max-stable
distribution with Smith dependence and fitted the conditioned
dependence model using:\ $i)$ the canonical family (\ref{eq:HT_class})
and $ii)$ the model implied by the norming functions (\ref{eq:a_hr}),\
treating the functions~(\ref{eq:a_hr}) as a parametric model for the
growth of $Y$ given large $X$.\ Our comparisons are based on the
differences between the conditional quantile estimates of $Y\mid X=x$
from the two models.
\begin{figure}[htbp!]
  \centering 
  \includegraphics[scale=.98,trim=30 20 30 30]{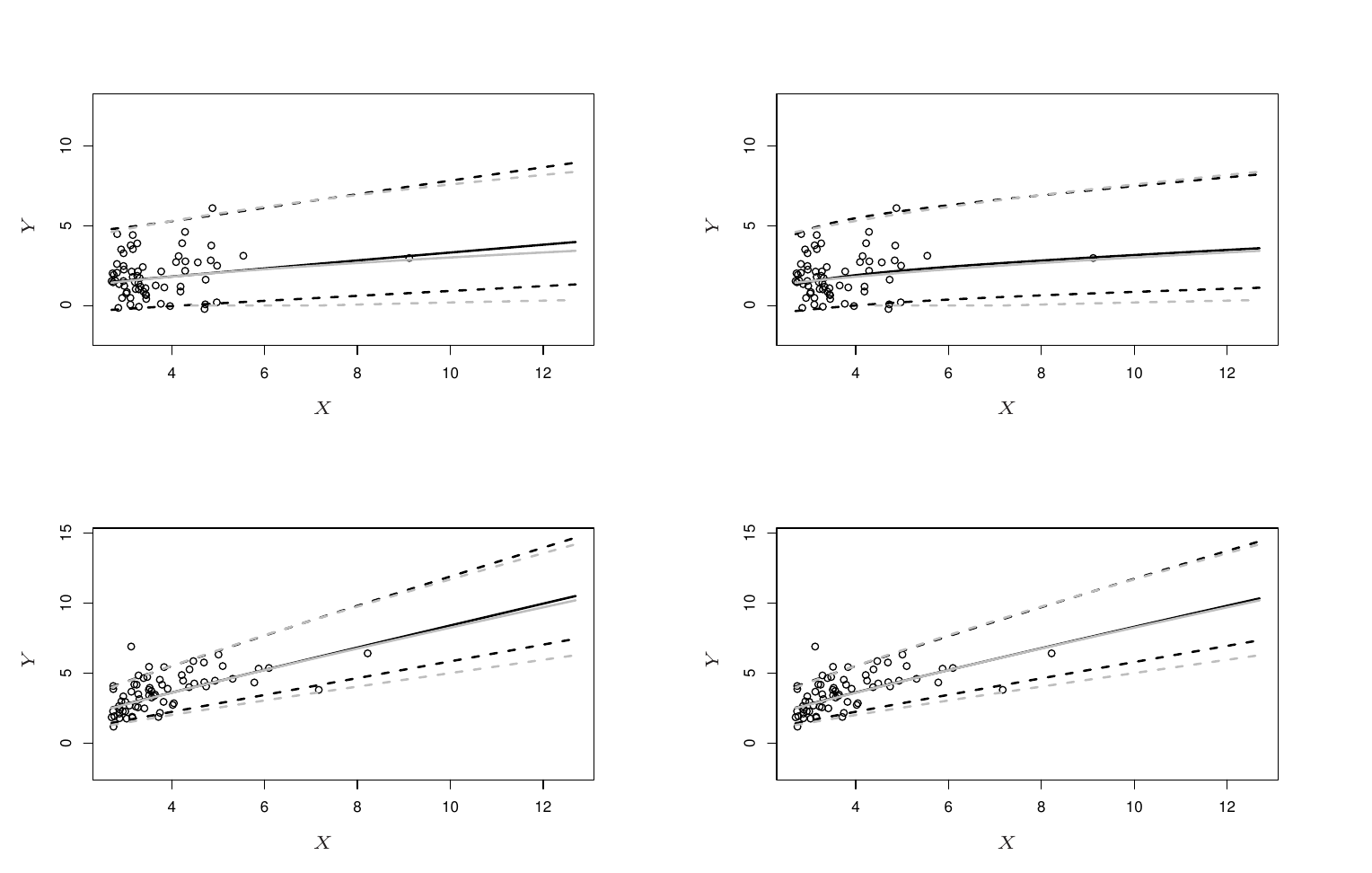}
  \caption{Conditional exceedances above the 0.935 unit exponential
    quantile from a simulated sample of size $10^3$ from the inverted
    max-stable distribution with Smith dependence and exponential
    margins for $\lambda=1.3$ (top) and $\lambda=0.3$ (bottom).\ The
    black lines correspond to averaged estimates from 100 simulations
    of the 0.025,\ 0.5 and 0.975 conditional quantiles of $Y\mid
    X=x$,\ $x>F_X^{-1}(0.935)$,\ using the Heffernan--Tawn model
    (\ref{eq:HT_class}) (left) and the model constructed by the
    theoretical functions in expression~(\ref{eq:a_hr}) (right).\ Grey
    lines correspond to the theoretical 0.025,\ 0.5 and 0.975
    conditional quantiles.}
\label{fig:sims}
\end{figure}
For both models, similar to \cite{hefftawn04}, we used for the
limiting law $G$ in expression~(\ref{eq:cond_rep}) with the false
working assumption of a normal distribution with mean and variance
parameters.\ We considered two values for the dependence parameter,\
i.e.,\ $\lambda=1.3$ (weak dependence) and $\lambda=0.3$ (strong
dependence).\ For each $\lambda$,\ $10^2$ samples of size $10^3$ were
generated from the inverted max-stable distribution with Smith
dependence and the 0.025,\ 0.5 and 0.975 conditional quantile
estimates of $Y\mid X=x$,\ for $x>F_X^{-1}(0.935)$,\ were computed
from the two model fits,\ i.e.,\ model~(\ref{eq:HT_class}) and the
model defined by expression (\ref{eq:a_hr}).\ The conditional quantile
estimates are of the form $\hat{a}(x) + \hat{b}(x) \hat{z}_p$,\ where
$\hat{a}(x)$,\ $\hat{b}(x)$ are maximum likelihood estimates and
$\hat{z}_p$ is the $p$-th empirical quantile of $\hat{Z}=\{Y -
\hat{a}(x)\}/\hat{b}(x)$,\ for large $x$.\ Figure~\ref{fig:sims} shows
the averaged estimates of the conditional quantiles along with the
theoretical conditional quantiles.\ Both models estimate the true
conditional quantiles well and their behaviour is almost
indistinguishable.\ This shows that the canonical model is flexible
enough to approximate the conditional distribution of the inverted
max-stable distribution with Smith dependence.

\subsection{$\Gamma$-variation at lower tail of spectral density}
\label{sec:new_model}
Having identified a new form for the tail of the spectral density for
the Smith max-stable model,\ we consider in this section the
log-conditional survivor function of $Y\mid X=x$,\ under the
assumption
\begin{equation}
  h(w) \sim g(w-w_{\ell}), \quad \text{as $w \rightarrow w_{\ell}$},
\label{eq:spec_gamma_varying}
\end{equation}
where $g(w)\in\Gamma_f(0^+)$.\ Similar to Section~\ref{sec:reg_var},\
we consider the assumption of possible mass at the lower end point
$w_{\ell}$.\ Our findings are based on the assumptions of a differentiable
spectral density $h$ and Lemma~\ref{le:gamma}.
\begin{lemma}
  \label{le:gamma} Let $g:\mbbR^+ \rightarrow \mbbR^+\in\Gamma_f(0^+)$
  and $U(w)\in R_{\nu}(0^+)$,\ $\nu\in \mbbR$.\ Assume further that
  there exists an $\epsilon>0$ such that $U$ and $g$ are
  $C^\infty(0,\epsilon)$ functions with
  \begin{equation}
    \lim_{w \rightarrow 0^+}
    \frac{g^{\prime}(w)} {g(w)^2} \int_{0}^{w}g(s)~\mathrm{d}s
    \label{eq:2ndlimit}
  \end{equation}
existing. Then\newline
  \noindent $(i)$~~~$U(w) g(w)\in\Gamma_f(0^+)$.\newline  
  \noindent $(ii)$~~Define $f(w) = g(w)/g'(w)$,\ $w>0$.\ Then,\ $f$ is an
  auxiliary function for $g$.\newline
  \noindent $(iii)$ For $f(w)$ as in $(ii)$,
  \begin{IEEEeqnarray}{rCl}
    \left(\int_{0}^w U(s) g(s)
      \mathrm{d}s\right)\big/\left\{U(w)f(w)g(w)\right\} &=& 1 -
    \frac{(U
      f)'(w)}{U(w)}  + \frac{\left\{ (Uf)'f\right\}'(w)}{U(w)} - \cdots    \label{eq:expansion}\\\nonumber\\
    &=& 1 + o(1), \quad \text{as $w \rightarrow
      0^+$} \label{eq:little_o}.
    \end{IEEEeqnarray}
\end{lemma}
\noindent Proposition~\ref{prop:gamma_variation} gives the asymptotic
form of the log-conditional survivor function.
\begin{proposition}
  \label{prop:gamma_variation}
  Under the conditions of Lemma~\ref{le:derivative} and for $h(w)$ as
  in expression~(\ref{eq:spec_gamma_varying}), for $y/(x+y)\downarrow
  w_{\ell}\in [0,1/2)$ for $y=y(x)$ function of $x$, as $x \rightarrow
  \infty$,\ we obtain that $\log\Pr\left(Y>y \mid X=x\right)$ is
  asymptotically equivalent to

\noindent $(i)$ for $H(\{w_{\ell}\})=0$,
\begin{IEEEeqnarray}{rCl}
 && - \left(x+y\right) f^2\left(\frac{y}{x+y}-w_{\ell}\right)
    h\left(\frac{y}{x+y}\right),
\label{eq:nomass_cond_distr2}
\end{IEEEeqnarray}
 
\noindent $(ii)$ for $H(\{w_{\ell}\}) > 0$, 
\begin{IEEEeqnarray}{rCl}
  && \log\left\{1 - w_{\ell} H(\{w_{\ell}\})\right\} - (x+y)\left(\frac{y}{x+y} -
    w_{\ell}
  \right) H(\{w_{\ell}\}).
\label{eq:mass_cond_distr2}
\end{IEEEeqnarray}
\end{proposition}

As an example,\ we explore a new class of spectral densities that are
more flexible than the spectral density~(\ref{eq:h_hr}) of the
inverted max-stable distribution with Smith dependence.\
Specifically,\ consider for $\gamma > 0$,\ $\delta\in \mathbb{R}$ and
$\kappa>0$,\ the $\Gamma_f(0^+)$-varying function
  \begin{IEEEeqnarray}{rCl}
    &h(w)\sim w^\delta \exp\left(-\kappa w^{-\gamma}\right)& \quad
    \text{as $w\rightarrow 0^+$},
    \label{eq:h_new}
  \end{IEEEeqnarray}
  with auxiliary function
  \begin{equation}
    f(w) = \left(\kappa \gamma\right)^{-1} w^{1+\gamma}.
    \label{eq:auxiliary_new}
  \end{equation}
  Proposition~\ref{prop:new_model} gives the normalising functions and
  limiting conditional distribution for this example.
\begin{proposition} 
\label{prop:new_model}
Assume that $(X,Y)$ follows the inverted max-stable
distribution~(\ref{eq:inv_max_stable}) with spectral measure $H$ with
$w_{\ell}=0$ and $H(\{0\})=0$,\ and spectral density $h$ that satisfies
expression~(\ref{eq:h_new}).\ Limit~\eqref{eq:cond_rep} holds with, for $x>0$ and $z \in \mathbb{R}$
\begin{equation}
  a(x) = x \kappa^{1/\gamma}\left(\log x\right)^{-1/\gamma}\left[ 1
    + \gamma^{-2} \left\{\delta + 2 (1+\gamma)\right\} \frac{\log \log
        x}{\log x}\right],\, b(x) = x \left(\log x\right)^{-1-1/\gamma},
  \label{eq:a_new}
\end{equation}
and
\begin{equation}
  G(z) = 1 -
  \exp\left[-\left\{\kappa^{(\delta+2)/\gamma}/\gamma^2\right\}\exp\left(\gamma
      \kappa^{-1/\gamma}z\right)\right] 
  .
\label{eq:new_models_limit}
\end{equation}
\end{proposition}
\begin{remark}
  Similarly with the inverted max-stable with Smith dependence,\ the
  limiting conditional distribution~\eqref{eq:new_models_limit} is of
  reverted Gumbel type and the norming functions~(\ref{eq:a_new}) do
  not belong to the Heffernan-Tawn parametric family.\ The rate of
  convergence to the limit is order $\log\log u/\log u$.
\end{remark}
\begin{figure}[htpb!]
  \centering
  \includegraphics[scale=1]{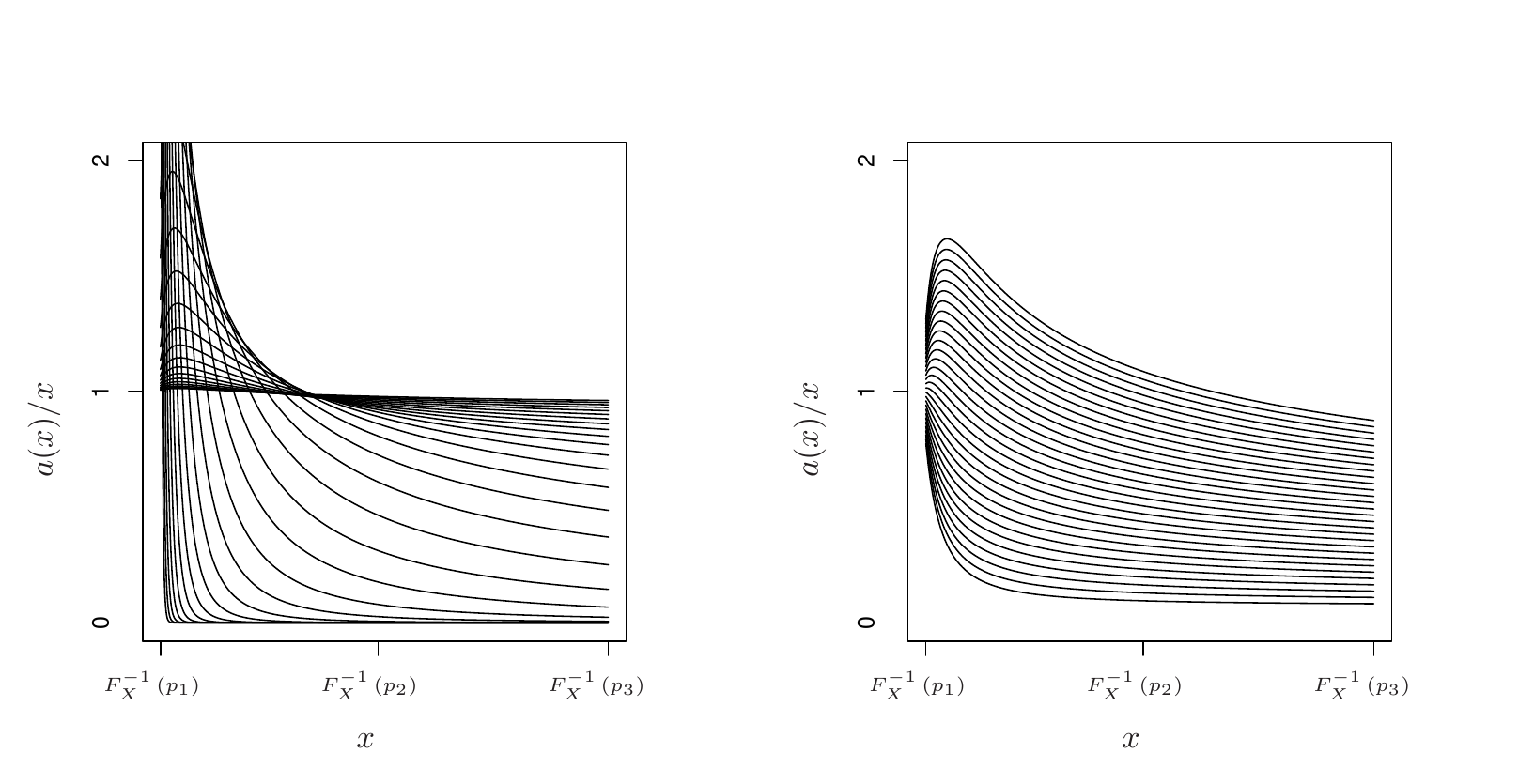}
  \includegraphics[scale=1]{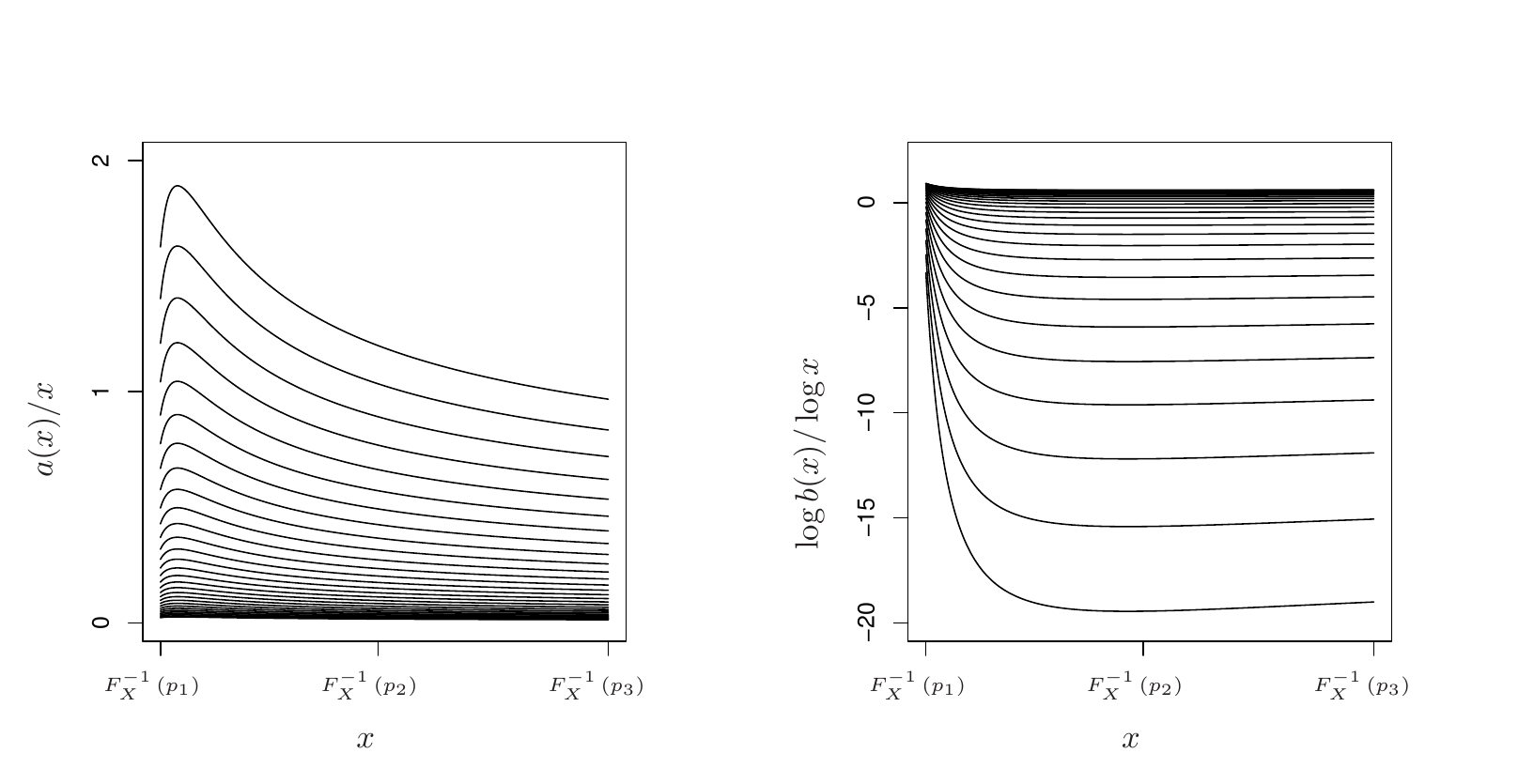}\
  \caption{Plots of $a(x)/x$ for different values of the parameters
    $\gamma$ (top left),\ $\kappa$ (top right) and $\delta$ (bottom
    left) and $\log b(x)/ \log x$ (bottom right),\
    $x>F_X^{-1}(0.87)$. The inverse of the unit exponential
    distribution function is $F_X^{\leftarrow}(p)=-\log(1-p)$,\ for
    $p\in(0,1)$.\ The values of $p_1$,\ $p_2$ and $p_3$ are 0.95,\
    $1-10^{-7}$ and $1-10^{-13}$,\ respectively.}
\label{fig:norm_functions2}
\end{figure}
\pagebreak
Figure~\ref{fig:norm_functions2} shows the graphs of
functions~$a(x)/x$ and $\log b(x)/\log x$,\ for several values of the
parameters $\gamma$,\ $\delta$,\ $\kappa$,\ and large $x$.\ Large
values of $\gamma$ correspond to strong dependence between $X$ and $Y$
so that $a(x)/x$ is nearly constant and equal to $1$ for large $x$.\
Small values of $\gamma$ correspond to independence with sharp
decrease of $a(x)/x$ to $0$ as $x$ increases.\ Intermediate values of
$\gamma$ correspond to mild-moderate dependence of $X$ and $Y$ with
$a(x)/x$ having a turning point and decaying as $x$ increases.\
Parameters $\kappa$ and $\delta$,\ seem to have similar effect with
larger values corresponding to increasing dependence.\ Last,\ $\log
b(x)/\log x$ is approximately constant with large $x$.\ Comparing with
the canonical family~(\ref{eq:HT_class}),\ the degree of approximation
of $a(x)$ and $b(x)$ by $\alpha x$ and $x^\beta$ is not especially
good for $a(x)$.\ To our knowledge,\ there is no current parametric
model for the max-stable process with spectral density satisfying
expression~(\ref{eq:h_new}).
\section{Spatial extensions}
\label{sec:spatial}
\cite{engeetal2014} have found considerable insights can be gained by
studying max-stable processes conditioned on the process being extreme
at a fixed point, $ s^* \in \mathbb{R}^2$ say.\ For $\{X(s); s\in
\mathbb{R}^2\}$ an inverted max-stable processes with unit exponential
margins a similar approach requires the existence of location and
scaling norming functions $a_s:\mathbb{R}_{+} \rightarrow \mathbb{R}$
and $b_s:\mathbb{R}_{+} \rightarrow \mathbb{R}_{+}$, for all $s\in
\mathbb{R}^2$, such that, for any $x>0$ and $\{z_s \in\mbbR, s\in
\mathbb{R}^2\}$,
\begin{equation}
  \lim_{u\rightarrow \infty}
  \Pr\left[X(s^*)-u>x,\frac{X(s)-a_s\{X(s^*)\}}{{b_s\{X(s^*)\}} } <z_s, s\in \mathbb{R}^2\backslash\{ s^*\}
    \mid
    X(s^*)>u\right]  = \exp(-x)G( \boldsymbol{z}),
  \label{eq:cond_rep_spat}
\end{equation}
where $G$ is an infinite-dimensional joint distribution function that
is non-degenerate in all univariate margins and $\boldsymbol{z}=\{z_s,
s\in \mathbb{R}^2\backslash\{ s^*\}\}$.

From our working in this paper we have derived the forms of $(a_s,
b_s)$ and $G_s$, the marginal of $G$, as these features are determined
uniquely by the bivariate joint distributions of the inverted
max-stable process.\ Thus all that remains to fully characterize
limit~\eqref{eq:cond_rep_spat} is to derive the infinite-dimensional
dependence structure of $G$. 

Here a complication arises as for most max-stable process models only
the bivariate marginal distributions are available in closed form
\citep{davpadrib}. Maybe future work can explore approaches that
\cite{engeetal2014} have used to get around this issue for max-stable
process.\ However,\ for the moment we restrict ourselves to consider
the one well-known max-stable process which has closed form trivariate
marginal distributions, the Smith process \citep[see][]{gentetal11}.\ In
this case we find that the associated bivariate limit
$G_{s_1,s_2}(z_1,z_2)= G_{s_1}(z_1) G_{s_2}(z_2)$, thus $G$
factorises; corresponding to asymptotic conditional independence in
the terminology of \cite{hefftawn04}.\ It follows that the
infinite-dimensional $G$ will also factorise.\ Thus in this case all
the dependence structure of the inverted max-stable processes is
absorbed in the location-scale functions. Identifying which classes of
max-stable process possess this asymptotic conditional independence
for $G$ is an interesting line of future research.


\section{Appendix}
\label{sec:appendix}


\subsection{Proof of Lemma~\ref{le:derivative}}

For $s/(s+t)\in (w_{\ell},w_{u})$, the partial derivative of the exponent
measure is equal to
\begin{IEEEeqnarray*}{rCl}
  \frac{\partial V (s,t)}{\partial s}& = & \frac{\partial }{\partial s} \left[ \int_{\frac{s}{s+t}}^{w_{u}}
    (w/s)\mathrm{d}H(w) + \int_{w_{\ell}}^{\frac{s}{s+t}}
    \{(1-w)/t\}\mathrm{d}H(w)\right] \\\\
  & = & \frac{\partial }{\partial s} \left[ \int_{\left[\frac{s}{s+t},w_{u}\right)}
    (w/s)h(w)\mathrm{d}w + \int_{\left(w_{\ell},\frac{s}{s+t}\right] }
    \{(1-w)/t\} h(w)\mathrm{d}w \right] -  (w_{u}/s^2) H(\{w_{u}\})\\\\
  & = & -\int_{\left[\frac{s}{s+t},w_{u}\right) } (w/s^2)h(w)\mathrm{d}w -
  (w_{u}/s^2) H(\{w_{u}\}),
\end{IEEEeqnarray*}
which, under the assumption of $y/(x+y) \downarrow w_{\ell}$, as
$x\rightarrow \infty$, yields 
\begin{equation}
  V_1(1,x/y) \rightarrow -\int_{w_{\ell}}^{ w_{u}} wh(w)\mathrm{d}w - 
  w_{u} H(\{w_{u}\}),
\label{eq:limit_deriv_appendix}
\end{equation}
as $x\rightarrow \infty$.\ Using the moment
constraint~(\ref{eq:max_stable}) we have that
\[
\int_{w_{\ell}}^{w_{u}} w h(w) \mathrm{d}w =
\int_{\left(w_{\ell},w_{u}\right)} w h(w) \mathrm{d}w =1 -
w_{\ell} H(\{w_{\ell}\}) - w_{u} H(\{w_{u}\}),
\]
which yields equation~(\ref{eq:partial_limit}), after combining with
equation~(\ref{eq:limit_deriv_appendix}).

\subsection{Proof of Proposition~\ref{prop:log_survivor_cases}}
Working similarly to the proof of Lemma~\ref{le:derivative}, we get,
after combining equations~(\ref{eq:approx_surv})
and~(\ref{eq:partial_limit}), that for $c(x,y)=y/\{x+y\}$,
$d(x,y)=\{x+y\}/\{(1-w_{\ell})y - w_{\ell} x\}$, the log-conditional survivor,
$\log\left\{\Pr\left(Y>y \mid X=x\right)\right\}$, is equal to\newline
\noindent $(i)$ for $H(\{w_{\ell}\}) = 0$, 
\begin{IEEEeqnarray*}{rCl}
  && \left(x+y\right)\int_{w_{\ell}}^{c(x,y)} w \mathcal{L}(w-w_{\ell})
  (w-w_{\ell})^{t}~\mathrm{d}w -
  y\int_{w_{\ell}}^{c(x,y)} \mathcal{L}(w-w_{\ell})(w-w_{\ell})^t ~\mathrm{d}w \\\\
  &=& \left(x+y\right)\int_{d(x,y)}^{\infty} \mathcal{L}(1/s)
  s^{-(t+3)}~\mathrm{d}s + \left\{w_{\ell} x - (1-w_{\ell})
    y\right\}\int_{d(x,y)}^{\infty} \mathcal{L}(1/s)s^{-(t+2)}
  ~\mathrm{d}s.
\end{IEEEeqnarray*}
For $t>-1$, $d(x,y)\rightarrow \infty$, $c(x,y)\rightarrow w_{\ell}$, as $x
\rightarrow \infty$, we have, from Karamata's theorem
\citep[pg.~17]{resn87}, that the last expression is asymptotically equivalent to
\begin{IEEEeqnarray*}{rCl}
  && \frac{\left(x+y\right)}{(t+2)} \left\{d(x,y)\right\}^{-(t+2)}
  \mathcal{L}\{1/d(x,y)\} + \frac{\left\{w_{\ell} x -
      (1-w_{\ell})y\right\}}{(t+1)} \left\{d(x,y)\right\}^{-(t+1)}
  \mathcal{L}\{1/d(x,y)\},
\end{IEEEeqnarray*}
as $x\rightarrow \infty$, which simplifies to
expresion~(\ref{eq:nomass_cond_distr}).\newline
\noindent $(ii)$ for $H(\{w_{\ell}\}) > 0$,
\[\log\left\{1 - w_{\ell} H(\{w_{\ell}\})\right\} + \left\{w_{\ell} x - (1-w_{\ell})
    y\right\} H(\{w_{\ell}\})\].
\subsection{Proof of Proposition~\ref{prop:HR}}
Let $\phi$ be the probability density function of the standard normal
distribution.\ Assuming $y\rightarrow \infty$ as $x\rightarrow \infty$
with $y/x \rightarrow 0$, we have from
expression~(\ref{eq:approx_surv}), Lemma~\ref{le:derivative} and
Mill's ratio, that for large $x$, the log-conditional survivor, $\log
\Pr\left(Y>y \mid X=x\right)$, is approximately equal to
\begin{IEEEeqnarray}{rCl}
  && x - x\left[1 - \frac{\phi\left\{\frac{\lambda}{2}+
      \frac{1}{\lambda}\log(x/y)\right\}}{\frac{\lambda}{2}+
    \frac{1}{\lambda}\log(x/y)}\right] + y
  \frac{\phi\left\{\frac{\lambda}{2}-
      \frac{1}{\lambda}\log(x/y)\right\}}{\frac{\lambda}{2}-
    \frac{1}{\lambda}\log(x/y)}
  \nonumber \\\nonumber\\
  &=& x \frac{\phi\left\{\frac{\lambda}{2}+
      \frac{1}{\lambda}\log(x/y)\right\}}{\frac{\lambda}{2}+
    \frac{1}{\lambda}\log(x/y)} \left[1 +
    \frac{y}{x}\frac{\phi\left\{\frac{\lambda}{2}-
        \frac{1}{\lambda}\log(x/y)\right\}}{\phi\left\{\frac{\lambda}{2}+
        \frac{1}{\lambda}\log(x/y)\right\}}
    \frac{\left\{\frac{\lambda}{2}+
        \frac{1}{\lambda}\log(x/y)\right\}} {\left\{\frac{\lambda}{2}-
        \frac{1}{\lambda}\log(x/y)\right\}}\right]
  \nonumber \\\nonumber\\
  &=& x \frac{\phi\left\{\frac{\lambda}{2}+
      \frac{1}{\lambda}\log(x/y)\right\}}{\frac{\lambda}{2}+
    \frac{1}{\lambda}\log(x/y)} \left[1 +
    \frac{\left\{\frac{\lambda}{2}+
        \frac{1}{\lambda}\log(x/y)\right\}} {\left\{\frac{\lambda}{2}-
        \frac{1}{\lambda}\log(x/y)\right\}}\right]  \nonumber \\\nonumber\\
  &\doteq& -c \left(x y\right)^{1/2}\frac{\phi\left\{
      \frac{1}{\lambda}\log(y/x)\right\}}{
    \left\{\frac{1}{\lambda}\log(y/x)\right\}^2}\left[1 +
    O\left\{\left(\log x\right)^{-1}\right\}\right],
\label{eq:proof_inv_hr_limit}
\end{IEEEeqnarray}
where $ c= \lambda \exp\left(-\lambda^2/8\right) $.\ Now, let $z\in
\mbbR$ and $y=a(x)+b(x)z$, where $a(x)$ and $b(x)$ are given by
equations~(\ref{eq:a_hr}).\ We have, as $x\rightarrow \infty$,
\begin{IEEEeqnarray}{rCl}
  (xy)^{1/2}&=& x \exp\left\{\lambda^2/4 -\frac{\lambda}{\sqrt{2}}
    (\log x)^{1/2}\right\} \left[ 1 + O\left\{\frac{\log \log
        x}{\left(\log x\right)^{1/2}}\right\}\right],
  \label{eq:as_xy}\\\nonumber\\
  \left\{\frac{1}{\lambda}\log(x/y)\right\}^2&=& 2 \log x \left[1 + O \left\{\left(\log x\right)^{-1/2}\right\}\right],
  \label{eq:as_logsq}\\\nonumber\\
  \phi\left\{\frac{1}{\lambda}\log(x/y)\right\} &=&
  (2\pi)^{-1/2}\exp\left[-\log x -\frac{\lambda^2}{8} +
    \frac{\lambda}{\sqrt{2}}(\log x)^{1/2}+\log\log x + \frac{\sqrt{2}
      z}{\lambda}\right]\nonumber\\  \label{eq:as_phi}\\
  &&\times\left[1 + O\left\{\frac{(\log \log x)^2}{\log x}\right\}\right].\nonumber
\end{IEEEeqnarray}
Combining
equations~(\ref{eq:proof_inv_hr_limit}),~(\ref{eq:as_xy}),~(\ref{eq:as_logsq})
and~(\ref{eq:as_phi}) we get
\[
\Pr\left\{Y<a(x) + b(x)z \mid X=x\right\} = 1 - \exp\left[ -
  \frac{\lambda}{(8\pi)^{1/2}}\exp\left\{\sqrt{2}z/\lambda\right\}\right]
+ O\left\{\frac{\log \log x}{\left(\log x\right)^{1/2}}\right\}.
\]
Last, direct application of statement~(\ref{eq:reg_cond}) yields the
result of Proposition~\ref{prop:HR}.

\subsection{Proof of Lemma~\ref{le:gamma}}
$(i)$ First, for any auxiliary function $f$, we have that
$\lim_{t\rightarrow 0^+}f(t)/t = 0$ \citep[see][Lemma 1.5.1]{haan70}.\
Next, for $U(w)=w^\nu\mathcal{L}(w)$, where $\mathcal{L}\in R_0(0^+)$
and $\nu \in \mbbR$,
\begin{IEEEeqnarray*}{rCl}
  \lim_{t\rightarrow 0^+} \frac{U\left\{ t + w f(t) \right\}g(\left\{
      t + w f(t) \right\})}{U(t)g(t)} &=& \lim_{t\rightarrow
    0^+}\frac{\mathcal{L}\left[ t\left\{ 1 + w f(t)/t\right\} \right]}{\mathcal{L}(t)}
  \frac{g\left\{ t + w f(t)  \right\}} {g(t)} \left\{ 1 + w f(t)/t\right\},\\\\
  &\rightarrow & \exp(w), \quad \text{$w > 0$}.
\end{IEEEeqnarray*}
\noindent $(ii)$ Theorem 1.5.2 in \cite{haan70} asserts that any
function $f_0$ that satisfies
\begin{equation}
f_0(w)\sim\left(\int_{0}^wg(s)ds\right)/g(w)\quad \text{as $w\rightarrow
  0^+$},
\label{eq:de_haan}
\end{equation}
is an auxiliary function for $g$.\ Additionally, Theorem 1.5.4 in
\cite{haan70} states that if $g\in \Gamma_f(0^+)$, then, as $w
\rightarrow 0^+$
\begin{equation}
\int_{0}^w \{g(s)\}^2~\mathrm{d}s \sim \frac{1}{2} g(w)\int_{0}^w
g(s)~\mathrm{d}s.
\label{eq:thm154}
\end{equation}
Given that limit~(\ref{eq:2ndlimit}) exists it follows by
l'H\^{o}pital's rule that
\begin{equation}
\lim_{w\rightarrow 0^+}
\frac{1}{2} g(w)\int_{0}^w g(s)~\mathrm{d}s \Big/ 
\int_{0}^w \{g(s)\}^2~\mathrm{d}s
=
\lim_{w\rightarrow 0^+}
\left(\frac{1}{2} g(w)^2 +\frac{1}{2}g^{\prime}(w) \int_{0}^w g(s)~\mathrm{d}s\right)\Big /
g(w)^2,
\label{eq:3rdlimit}
\end{equation}
since the right hand side limit in~(\ref{eq:3rdlimit}) exists due to
limit~(\ref{eq:2ndlimit}) existing. However, by
property~(\ref{eq:thm154}), we have the left hand side limit
in~(\ref{eq:3rdlimit}) is equal to $1$ and so it follows from
equality~(\ref{eq:3rdlimit}) that as $w\rightarrow 0^+$,
\begin{equation}
  \int_{0}^w g(s)~\mathrm{d}s \sim \{g(w)\}^2/g'(w).
  \label{eq:eq:thm154diff}
\end{equation}
Combining expressions~(\ref{eq:de_haan}) and~(\ref{eq:eq:thm154diff})
we obtain $f_0(w) \sim f(w)$, as $w\rightarrow 0^+$, where
$f(w)=g(w)/g'(w)$.\ Hence,\ for $w>0$,\ the function $f(w)$ is up to
asymptotic equivalence equal to $f_0$.\newline

\noindent $(iii)$ Define $f(w)=g(w)/g'(w)$, $w>0$.\ We have
\begin{IEEEeqnarray*}{rCl}
  \int_{0}^{w} U(s)g(s)~\mathrm{d}s&=& \int_{0}^{w} U(s) f(s)
  g'(s)~\mathrm{d}s\\\\
  &=& U(w) f(w) g(w) - \int_{0}^{w} (Uf)'(s)
  f(s) g'(s)~\mathrm{d}s \\\\
  &=&U(w) f(w) g(w) - (Uf)'(w) f(w) g'(w) + \int_{0}^{w} \left\{(Uf)'f\right\}'(s) f(s)
  g'(s)~\mathrm{d}s,
\end{IEEEeqnarray*}
which gives expression~(\ref{eq:expansion}), after continuation of
integration by parts and division by $U(w)f(w)g(w)$.\ Last,
expression~(\ref{eq:little_o}) follows from de Haan's
theorem~(\ref{eq:de_haan}) and case $(i)$ of Lemma~\ref{le:gamma}.

\subsection{Proof of Proposition~\ref{prop:gamma_variation}}
We only consider case $(i)$ and note that case $(ii)$ is identical to
Proposition~\ref{prop:log_survivor_cases}.\ Define $c(x,y) = y/\left(x
  + y\right)$ and $l(x,y)=c(x,y)-w_{\ell}$.\ For $c(x,y)\rightarrow w_{\ell}$,
$l(x,y)=c(x,y)-w_{\ell} \rightarrow 0$, as $x \rightarrow \infty$, and with
$T$ as in expression~(\ref{eq:inv_max_stable}), we have that $ \log
\Pr\left(Y>y\mid X=x\right)$ is equal to
\begin{IEEEeqnarray}{rCl}
  && \left(x + y\right) \int_{0}^{l(x,y)}s g(s)~\mathrm{d}s +
  \left\{w_{\ell} x - (1-w_{\ell})y\right\} \int_{0}^{l(x,y)}
  g(s)~\mathrm{d}s.\label{eq:gamma_pr_min_A_min_B}
\end{IEEEeqnarray}
Using the asymptotic expansion~(\ref{eq:expansion}), up to first
order, we have that the two integrals in
expression~(\ref{eq:gamma_pr_min_A_min_B}) are, as $x\rightarrow
\infty$, asymptotically equivalent to
\begin{IEEEeqnarray}{lCl}
  \int_{0}^{l(x,y)}s g(s)~\mathrm{d}s &\doteq&
  l(x,y)f\left\{l(x,y)\right\} g\left\{l(x,y)\right\}
  \left[1 - \frac{f\left\{l(x,y)\right\} - l(x,y)f'\left\{l(x,y)\right\}}{l(x,y)}\right]\label{eq:int_A},\\\nonumber\\
  \int_{0}^{l(x,y)} g(s)~\mathrm{d}s &\doteq& 
  f\left\{l(x,y)\right\} g\left\{l(x,y)\right\}
  \left[1 - f'\left\{l(x,y)\right\}\right],
  \label{eq:int_B}
\end{IEEEeqnarray}
for $U(s)=s \in R_1(0^+)$ and $U(s)=1\in R_0(0^+)$, respectively.
Combining expressions~(\ref{eq:gamma_pr_min_A_min_B}),
(\ref{eq:int_A}) and (\ref{eq:int_B}) we get
that~(\ref{eq:gamma_pr_min_A_min_B}) is asymptotically equivalent to
\begin{IEEEeqnarray*}{rCl}
  && f\left\{l(x,y)\right\} g\left\{l(x,y)\right\} \left(-\left\{w_{\ell} x
      - (1-w_{\ell})y\right\}\left[f'\left\{l(x,y)\right\} -
      \frac{f\left\{l(x,y)\right\}-l(x,y)f'\left\{l(x,y)\right\}}
      {l(x,y)}\right] \right)\\\\
  &=& f\left\{l(x,y)\right\} g\left\{l(x,y)\right\}\left[-\{w_{\ell} x -
    (1-w_{\ell})y\} \frac{(x + y) f\left\{l(x,y)\right\} }{w_{\ell} x -
      (1-w_{\ell})y}\right]
  \\\\
  &=& -\left(x+y\right) f^2\left(\frac{y}{x+y}-w_{\ell}\right)
  h\left(\frac{y}{x+y}\right),
\end{IEEEeqnarray*}
which completes the proof.

\subsection{Proof of Proposition~\ref{prop:new_model}}

Assuming $y\rightarrow \infty$ as $x\rightarrow \infty$ with $y/x
\rightarrow 0$, we have from expressions~(\ref{eq:approx_surv}),
(\ref{eq:nomass_cond_distr2}), (\ref{eq:h_new}),
(\ref{eq:auxiliary_new}) and Lemma~\ref{le:derivative}, that for large
$x$, the log-conditional survivor, $\log \Pr\left(Y>y \mid
  X=x\right)$, is approximately equal to
\begin{IEEEeqnarray}{rCl}
  && -\frac{x}{(\kappa \gamma)^2}
  \left(y/x\right)^{\delta+2(1+\gamma)}\exp\left\{-\kappa
    \left(y/x\right)^{-\gamma}\right\}.
\label{eq:gamma_limit}
\end{IEEEeqnarray}
Now, let $z\in \mbbR$ and $y=a(x)+b(x)z$, where $a(x)$ and $b(x)$ are
given by equations~(\ref{eq:a_new}).\ We have, as $x\rightarrow
\infty$,
\begin{IEEEeqnarray}{rCl}
  (y/x)^{\delta+2(1+\gamma)}&=& \left(\frac{\log
      x}{\kappa}\right)^{-\{\delta+2(1+\gamma)\}/\gamma}\left\{1+O\left(\frac{\log
        \log x}{\log x}\right)\right\},
  \label{eq:as_yovx}\\\nonumber\\
  \exp\left\{-\kappa \left(y/x\right)^{-\gamma}\right\}&=&x^{-1}
  \left(\frac{\log x}{\kappa}\right)^{\{\delta +
    2(1+\gamma)\}/\gamma}\exp\left( \gamma
    \kappa^{-1/\gamma}z\right)\left[1 + O\left\{\left(\frac{\log \log
          x}{\log x}\right)^2\right\}\right].
  \label{eq:expyovx}
\end{IEEEeqnarray}
Combining equations~(\ref{eq:gamma_limit}),~(\ref{eq:as_yovx})
and~(\ref{eq:expyovx}) we get
\[
\Pr\left\{Y<a(x) + b(x)z \mid X=x\right\} = 1 -
\exp\left[-\left\{\kappa^{(\delta+2)/\gamma}/\gamma^2\right\}\exp\left(\gamma
    \kappa^{-1/\gamma}z\right)\right] + O\left(\frac{\log \log
    x}{\log x}\right).
\]
Last, direct application of statement~(\ref{eq:reg_cond}) yields the
result of Proposition~\ref{prop:new_model}.


\subsection*{Acknowledgments}
Ioannis Papastathopoulos acknowledges funding from the SuSTaIn program
- Engineering and Physical Sciences Research Council grant
EP/D063485/1 - at the School of Mathematics of the University of
Bristol. We would like to thank Anthony C. Davison and Jennifer
L. Wadsworth for helpful discussions that improved the clarity of the
paper.

\end{document}